\documentclass[11pt,reqno]{amsart}
\usepackage{amssymb,amsmath,amsthm,amsfonts}

\newcommand{\T}{\mathbb{T}}
\newcommand{\N}{\mathbb{N}}

\newcommand{\C}{\mathbb{C}}

\newcommand{\D}{\mathbb{D}}
\newcommand{\Z}{\mathbb{Z}}
\newcommand{\cQ}{{\mathcal{Q}}}

\newcommand{\cS}{{\mathcal{S}}}

\newcommand{\aA}{{\bf\mathcal{A}}}
\newcommand{\cM}{{\mathcal{M}}}
\newcommand{\cun}{{\bf\mathcal{C}_1}}

\newcommand{\cn}{{\bf\mathcal{C}_{\it n}}}

\newtheorem*{teto*}{{\bf Th\'eor\`eme}}
\newtheorem*{cor}{\bf Corollaire}

\newtheorem{thm}{{\bf Th\'eor\`eme}}

\newtheorem{deff}[thm]{D\'efinition}

\newtheorem{lem}[thm]{Lemme}
\newtheorem{slem}[thm]{Sous--lemme}

\begin{document}

\title[Ideaux fermes d'algebres de Beurling]{
Ideaux fermes d'algebres  de Beurling   analytiques sur le bidisque}

\author[Bouya; El--Fallah; Kellay]{B. Bouya; O. El--Fallah \& K. Kellay}
\address{Brahim Bouya \& Omar El--Fallah. Departement de Mathematiques et d'Informatique. Universite Mohamed V, B.P. 1014. Rabat Maroc}
\curraddr{Brahim Bouya. Laboratoire Paul Painlev\'e, Universit\'e des Sciences et Technologies de Lille, B\^at. M2, 59655 Villeneuve d'Ascq Cedex, France}
\email{brahimbouya@gmail.com,  elfallah@fsr.ac.ma}
\address{Karim Kellay. LATP-CMI. Universite de Provence. 39 Rue F.Joliot Curie. 13453 Marseille cedex 13       France}
\email{kellay@cmi.univ-mrs.fr}
\maketitle

\subsection*{ \small Abstract}{\small
We study the closed ideal in  the Beurling algebras $\aA^{+}_{\alpha,\beta}$  of
 holomorphic function $f$ in the bidisc such that $\sum_{n,m\geq 0}
|\widehat{f}(n,m)|(1+n)^{\alpha}(1+m)^\beta<+\infty$.  We determine the function $f\in\aA^{+}_{\alpha,\beta}$  such that the ideals generated by $f$ coincide with the
ideal generated by their zeros set.}

\section{{{Introduction}}}
Soient $\D$ le disque unit\'e du plan complexe et $\aA({\D ^n})$
l'alg\`ebre du polydisque, l'alg\`ebre des fonctions continues sur
${\overline \D }^n$ et holomorphes sur $\D ^n$ munie de la norme
$$\| f\|_ {\infty}= \sup\{ |f(z)|:\ z\in \D ^n\}.$$
Lorsque $B$ est une alg\`ebre de Banach incluse dans $\aA({\D ^n})$, $f\in B$ et $E$ est un ensemble ferm\'e de ${ \overline \D}^n$; on notera par:
\begin{itemize}
\item $I_B(f)= { \overline{ fB}}$ l'id\'eal ferm\'e de $B$ engendr\'e par $f$.
\item $Z_f = \{ z\in { \overline \D }^N: \ f(z) =0\}$, l'ensemble des z\'eros de $f$.
\item $I_B(E) = \{ g\in B : \ g_{|E} = 0 \}$, l'id\'eal d'annulation de $B$ sur  $E$. 
\end{itemize}
Dans ce travail nous nous int\'eressons \`a la d\'etermination des fonctions $f \in B$ qui satisfont
$I_B(f)= I_B (Z_f)$ dans le cas o\`u $B$ est une alg\`ebre de Beurling analytique \`a poids polyn\^omial. Dans le cas de l'alg\`ebre du disque, $B= \aA({\D})$, le th\'eor\`eme de Beurling-Rudin donne une
 caract\'erisation compl\`ete
des id\'eaux ferm\'es de $\aA({\D})$. En particulier, si $f\in
\aA({\D})$ telle que $Z_f\subset \T$ alors $I_{\aA({\D})}(f) =
I_{\aA({\D})}(Z_f)$
 si et seulement si, $f$ est ext\'erieure~:
 $$f(z)=\exp\int_{0}^{2\pi}\frac{e^{i\theta}+z}{e^{i\theta}-z}\log|f(e^{i\theta})|\frac{d\theta}{2\pi}.$$
 Dans une s\'erie d'articles \cite{Hed1,Hed2,Hed3} motiv\'es par le prob\`eme de Levin
  \cite{Lev},  H. Hedenmalm s'est int\'eress\'e aux id\'eaux ferm\'es de certaines alg\`ebres de fonction en plusieurs variables. Il a obtenu dans le cas de l'alg\`ebre du
 bidisque les r\'esultats suivants~:
\begin{teto*}{\bf \cite{Hed1,Hed2}}
\begin{enumerate}
\item Soit $f \in \aA({\D ^2})$ telle que $Z_f\subset\{1\}\times\overline{\D}$, alors 
$I_{\aA({\D ^2})}(f)=I_{\aA({\D ^2})}(Z_f)$ si et seulement si les fonctions $f(\cdot,w)$
sont ext\'erieures pour tout $w\in\overline{\D}$ et la fonction
$f(1,\cdot)$ est soit identiquement nulle soit ext\'erieure.
\item   Soit $f \in \aA({\D ^2})$ telle que
$Z_f=\{1\}\times\overline{\D}\cup \overline{\D}\times\{1\},$
 si
$$|\log{|f(z,w)|}|=o(1/\inf(|1-z|,| 1-w|)),\quad
z\to 1\text{ ou }w\to 1 ,$$ alors $I_{\aA({\D ^2})}(f)=I_{\aA({\D ^2})}(\{1\}\times\overline{\D}\cup \overline{\D}\times\{1\}).$
\end{enumerate}
\end{teto*}
Pour $\alpha\geq 0$, consid\'erons maintenant les alg\`ebres de Beurling analytiques suivantes~:
$$ \aA^+_\alpha:=\{ f=\displaystyle \sum _{n\geq 0}a_nz^n:\ \|f\|_\alpha\displaystyle= \sum _{n\geq 0}|a_n|(1+n)^{\alpha} <\infty \}.$$
Dans \cite{Kah}, J.P. Kahane a montr\'e que si $f \in \aA^+ _\alpha$ telle que $Z_f =\{ 1\}$ alors
$I_{\aA^+_\alpha}(f) = I_{\aA^+_\alpha}(\{1\})$ si et seulement si $f$ est ext\'erieure.
 Notons que dans ce cas cette condition est \'equivalente \`a
$$
\displaystyle \lim _{|z|\to 1^-}(1-|z|)\log |f(z)| =0.
$$
Mentionnons \'egalement que la caract\'erisation des id\'eaux de $\aA^+_\alpha$ 
parait compliqu\'e, voir \`a ce sujet \cite{ESZ}.

Pour tout $\alpha, \beta \geq 0$, on consid\`ere les alg\`ebres de
 Beurling du bidisque suivantes~:
\begin{multline*}\aA^{+}_{\alpha,\beta}:=\big\{f (z,w)=\displaystyle \sum_{n,m\geq 0}a_{n,m}z^nw^m\in\aA({\D}^{2}) \text{ : }\\
\| f\|_{{\alpha,\beta}}:=\sum\limits_{n,m\in\N}
|a_{n,m}|(1+n)^{\alpha}(1+m)^\beta<+\infty\big\},
\end{multline*}
Dans ce travail nous \'etendons les r\'esultats obtenus par H.Hedenmalm (\cite{Hed1,Hed2,Hed3})
 aux alg\`ebres $\aA^{+}_{\alpha,\beta}$.
 Nous montrons les deux th\'eor\`emes suivants:

\begin{thm}\label{theo11} Soient $(\alpha,\beta)\in[0,1[\times]0,1[$ et $f\in\aA^{+}_{\alpha,\beta}$
 telle que $Z_f =\{1\}\times\overline{\D}$. Alors
$I_{\aA^{+}_{\alpha,\beta}}(f)=I_{\aA^{+}_{\alpha,\beta}}(\{1\}\times\overline{\D})$
si et seulement si les fonctions $f(\cdot,w)$
sont ext\'erieures pour tout $w\in\overline{\D}$.
\end{thm}
Comme cons\'equence du  Th\'eor\`eme \ref{theo11} nous obtenons 
\begin{cor}
Soient $(\alpha,\beta)\in[0,1[\times]0,1[$ et $f\in\aA^{+}_{\alpha,\beta}$
 telle que $Z_f =\{(1,1)\}$. Alors
$I_{\aA^{+}_{\alpha,\beta}}(f)=I_{\aA^{+}_{\alpha,\beta}}(\{(1,1)\})$ si et seulement si les fonctions $f(\cdot,1)$ et $f(1,\cdot)$
sont ext\'erieures.
\end{cor}

\begin{thm}\label{theo12}Soient $(\alpha,\beta)\in]0,1[\times]0,1[$ et $f\in
\aA^{+}_{\alpha,\beta}$ telle  que  $Z_f=\{1\}\times{{\overline
\D}}\cup \overline{{\D}} \times \{1\}.$ Si
$$|\log{|f(z,w)|}|=o(1/\inf(|1-z|,| 1-w|)),\qquad
z\to 1\text{ ou }w\to 1 ,$$ alors
$I_{\aA^{+}_{\alpha,\beta}}(f)=I_{\aA^{+}_{\alpha,\beta}}(\{1\}\times{\overline{{\D}}}\cup \overline{{\D}} \times \{1\})$.
\end{thm}
Notons que pour $\alpha \geq 1$ ou $\beta \geq 1$, des r\'esultats
analogues, faisant intervenir les d\'eriv\'ees partielles de $f$, peuvent \^etre obtenus de la m\^eme mani\`ere.

Le paragraphe 2 sera consacr\'e \`a l'\'etude de la notion de la
$\delta$-visibilit\'e qui jouera un r\^ole fondamental dans la
preuve du th\'eor\`eme 1 et du th\'eor\`eme 2.


\section{La $\delta$--visibilite}

Dans cette partie nous rappelons la notion de la
$\delta$--visibilit\'e introduite et \'etudi\'ee dans \cite{AD,EE,ENZ,EZ,N}.
 Soit $A$ une alg\`ebre de Banach commutative
unitaire. L'ensemble des caract\`eres de $A$ sera not\'e $\mathcal{M}_{A}.$
La transformation de Gelfand associ\'ee \`a $A$ est l'application:
\begin{eqnarray}\label{Guelfand}
\mathcal{G}_A:&A&\longrightarrow \mathcal{C}(\mathcal{M}_{A})\\
\nonumber &x&\longrightarrow \widehat{x}:\phi\in
\mathcal{M}_{A}\longrightarrow \widehat{x}(\phi)=\phi(x),
\end{eqnarray}
o\`u $\mathcal{C}(\mathcal{M}_{A})$ est l'alg\`ebre des fonctions
continues sur $\mathcal{M}_{A}.$\\
Soit $n\in\N.$ Pour tout $f=(f_{1},f_{2},...,f_{n})\in A^{n},$ on
pose
$$\left\{
\begin{array}{lll}
\displaystyle\|f\|&:=&\displaystyle\Big(\sum_{k=1}^{n}\|f_{k}\|^{2}\Big)^{1/2}\\
\displaystyle \delta_{f}&:=&\displaystyle\inf_{\phi\in\cM_{A}}\Big(\sum_{k=1}^{n}\mid
\widehat{f_{k}}(\phi)\mid^{2}\Big)^{1/2}
\end{array}
\right.
$$
\begin{deff}\label{def:4} Soit $0<\delta\leq1$. On dit que
le spectre de $A$ est $\delta$--$n$--visible s'il existe une constante
$\cn(\delta),$ qui d\'epend  de $\delta$ et de $n$, telle que
pour tout $f=(f_{1},f_{2},...,f_{n})\in A^{n}$, $\|f\|\leq 1$ v\'erifiant
$\delta_{f}\geq\delta$,
il existe $g=(g_{1},g_{2},...,g_{n})\in A^{n}$ tels que
 $$
 \left\{\begin{array}{l}
 \sum\limits_{k=1}^{n}f_{k}g_{k}=1,\\\\
 \|g\|\leq \cn(\delta).\end{array}\right.
 $$
\end{deff}

Posons
$$\cn(\delta,A):=\sup_{f}\Big\{\inf\big\{\|g\|\text{ :
}g=(g_1,\cdots,g_n)\in A^{n}\text{ et } \sum_{k=1}^{n}g_{k}f_{k}=1\big\}\Big\},$$ le
sup \'etant pris sur les $f=(f_1,\cdots,f_n)\in A^{n}$ telles que
$\delta_{f}\geq\delta$ et $\|f\|\leq 1$.
En particulier
$$\cun(\delta,A):=\sup\{\|f^{-1}\|\text{ : }\|f\|\leq1\text{ et }
|\widehat{f}(\phi)|\geq\delta\quad(\phi\in\cM_A)\}$$


\subsection{ La $\delta$--$1$--visibilit\'e pour les alg\`ebres de Beurling} 
Consid\'erons les alg\`ebres de Beurling \`a poids suivantes
  $$
\aA_\alpha =\Big\{f\in\mathcal{C}({\T}) \text{ : } \|
f\|_{\alpha}:=\sum\limits_{n\in\Z}
|\widehat{f}(n)|(1+|n|)^{\alpha}<+\infty\Big\} .$$ Munie du produit
ponctuel et de la norme $\|.\|_{\alpha}$, $ \aA_\alpha $ est une
alg\`ebre de Banach commutative unitaire. En identifiant son spectre
\`a $ \T$, la transformation de Gelfand devient~:
$\mathcal{G}_{\aA{_\alpha }}(f) = f$ pour tout $f\in \aA_\alpha $.
Il a \'et\'e d\'emontr\'e dans  \cite{EE, ENZ}, que le spectre de
$\aA_\alpha $ est $\delta$--$1$--visible
 pour tout $0<\delta \leq 1$. Plus pr\'ecis\'ement pour tout $\alpha \in ]0,1[$
 on a $\cun(\delta,\aA_\alpha)\leq c/\delta^{2+1/\alpha}$.\\
 Ce r\'esultat peut \^etre \'etendu aux alg\`ebres de Beurling
 en plusieurs variables. Soient $\alpha, \beta >0$
et soit $\aA_{\alpha,\beta}$ l'alg\`ebre de Beurling d\'efinie par:
 $$\aA_{\alpha,\beta}=\Big\{f\in\mathcal{C}({\T}^{2}) \text{ : } \|
f\|_{{\alpha,\beta}}:=\sum\limits_{n, m\in\Z^2}
|\widehat{f}(n,m)|(1+|n|)^{\alpha}(1+|m|)^\beta<+\infty\Big\} , $$
Munie du produit ponctuel et de la norme $\| .\|_{{\alpha,\beta}}$,
$\aA_{\alpha,\beta}$ est une alg\`ebre de Banach commutative et
unitaire dont le spectre peut \^etre identifi\'e \`a $\T ^2$. Nous
avons le r\'esultat suivant:

 \begin{lem}\label{symet2}Soient $\alpha>0$ et $\beta>0$. Posons
$\gamma=\inf\{1,\alpha,\beta\}.$ On a
$$
 \cun(\delta,\aA_{\alpha,\beta})\leq c\delta^{-(3+\frac{1}{\gamma}(1+\alpha+\beta))},\qquad 0<\delta<1,\\
$$  o\`u $c$
est une constante qui d\'epend seulement de $\alpha$ et de $\beta.$
\end{lem}

\begin{proof}
Soit $f\in\aA_{\alpha,\beta}$ telle que $|f|\geq \delta>0$ et telle que
$\|f\|_{\alpha,\beta}\leq1.$ Pour tout $0<\rho<1,$ on pose 
$$f_\rho(z,w):=\sum_{n,m\in\Z}\widehat{f}(n,m)\rho^{|n|+|m|}z^{n}w^{m},\qquad  \rho\leq|z|\leq1/\rho
\text{ et }\rho\leq|w|\leq1/\rho.$$
Posons $a_{n,m}=\widehat{f}(n,m)$ et $\gamma=\inf\{\alpha,\beta,1\}$. Pour $z,w\in\T$,  on a
\begin{multline*}
|f(z,w)-f_\rho(z,w)|=|\sum_{n,m\in\Z}a_{n,m}(1-\rho^{|n|+|m|})z^{n}w^{m}|\\
\leq
\sup_{n,m\in\Z}\frac{1-\rho^{|n|+|m|}}{(1+|n|)^{\alpha}(1+|m|)^{\beta}}
\leq 2(1-\rho)^{\gamma}.
\end{multline*}
Soit $\rho$ tel que $2(1-\rho)^{\gamma}=\delta/3$. Le rayon spectral
de $(f-f_\rho)f_{\rho}^{-1}$  est strictement inf\'erieur \`a $1$ et
\begin{equation}\label{del0}f^{-1}=\sum_{p\in\N}(f-f_\rho)^{p}f_{\rho}^{-p-1}.
\end{equation}
Soit $p\geq2,$ on a
\begin{multline}
  \|(f-f_\rho)^p\|_{\alpha,\beta}=\|\Big(\sum_{n,m\in\Z}a_{n,m}(1-\rho^{|n|+|m|})z^{n}w^{m}\Big)^p\|_{\alpha,\beta}\\
\leq
\sup_{n_i,m_j\in\Z}\frac{(1-\rho^{|n_1|+|m_1|})\ldots(1-\rho^{|n_p|+|m_p|})
}{(1+|n_1|)^\alpha(1+|m_1|)^\beta\ldots(1+|n_p|)^\alpha(1+|m_p|)^\beta}\times\\
\times(1+|n_1+\ldots+n_p|)^\alpha(1+|m_1+\ldots+m_p|)^\beta
\end{multline}
Puisque $(1+|n_1+..+n_p|)\leq p \displaystyle \max _{1\leq k\leq p}(1+|n_k|)$,
\begin{multline}
  \|(f-f_\rho)^p\|_{\alpha,\beta}\\\leq p^{\alpha+\beta}
\sup_{n_i,m_j\in\Z}\sup_{k=1}^{p}\sup_{l=1}^{p}\prod_{q=1}^{p}
\frac{(1-\rho^{|n_q|+|m_q|})}{(1+|n_q|)^\alpha(1+|m_q|)^\beta}(1+|n_k|)^\alpha(1+|m_l|)^\beta\\
\label{del2}\leq
p^{(\alpha+\beta)}(2(1-\rho)^\gamma)^{p-2}.
\end{multline}
Posons $\widehat{f_{\rho}^{-p-1}}(n,m)=b_{n,m},$   $p\geq
0$. L'in\'egalit\'e de H\"older et l'\'egalit\'e de Plancherel--Parseval, nous donne
\begin{multline}
   \|f_{\rho}^{-p-1}\|_{\alpha,\beta}\leq\\
\Big(\sum_{n,m\in\Z}\rho^{2(|n|+|m|)}(1+|n|)^{2\alpha}(1+|m|)^{2\beta}\Big)^{1/2}
\Big(\sum_{n\in\Z}\Big(\sum_{m\in\Z}|b_{n,m}|^2\rho^{-2|n|-2|m|}\Big)^{1/2}\\
\leq
\frac{c}{(1-\rho)^{\alpha+\beta+1}}\Big(\int_{\rho\T\cup\rho^{-1}\T}\int_{\rho\T\cup\rho^{-1}\T}
|f_{\rho}^{-p-1}(\xi,\zeta)|^2d\xi d\zeta\Big)^{1/2}\\
\leq\frac{c}{\delta^{p+1}(1-\rho)^{\alpha+\beta+1}}\label{del3}.
\end{multline}

De \eqref{del0}, \eqref{del2} et \eqref{del3} on en d\'eduit que
\begin{multline*}
\|f^{-1}\|_{\alpha,\beta}\leq\sum_{p\geq0}\|(f-f_\rho)^{p}\|_{\alpha,\beta}\|f^{-p-1}_\rho\|_{\alpha,\beta}
\\ \leq \|f^{-1}_\rho\|_{\alpha,\beta}+ \|f^{-2}_\rho\|_{\alpha,\beta} +\frac{c}{2(1-\rho)^{2\gamma+\alpha+\beta+1}}\sum_{p\geq2}
\frac{p^{(\alpha+\beta)}(2(1-\rho)^\gamma)^{p}}{\delta^{p+1}}\\
\leq
\frac{c}{\delta^{3}(1-\rho)^{\alpha+\beta+1}}=\frac{c}{\delta^{3+\frac{1}{\gamma}(\alpha+\beta+1)}}.
\end{multline*}
\end{proof}
\subsection{ La $\delta$--$2$--visibilit\'e de l'alg\`ebre $\aA^{+}_{\alpha,\beta}$}
 La $\delta$--$2$--visibilit\'e pour une alg\`ebre de Banach commutative et unitaire revient \`a la r\'esolution de l'identit\'e de
Bezout avec contr\^ole des normes des solutions. Le cas d'une variable \`a \'et\'e trait\'e dans \cite{EZ}.
Nous \'etendons ce r\'esultat au cas de plusieurs variables.
\begin{thm}\label{2-d-ab}
Soient
$f_1,f_2\in\aA^{+}_{\alpha,\beta}$ et $0<\delta$ telles que
  $$
\delta^2\leq  |f_1(z)|^2+|f_2(z)|^2\leq 1 , \qquad z\in\D^2.$$
Il existe
$h_1,h_2,\in \aA^{+}_{\alpha,\beta}$ telles que
$$\left\{
\begin{array}{l}f_1h_1+f_2h_2=1,\\ \\
\|h_1\|^{2}_{\alpha,\beta}+\|h_2\|^{2}_{\alpha,\beta }\leq \delta^{-c_{\alpha,\beta}}.\\
\end{array}
\right.
$$
 o\`u $c_{\alpha,\beta}$ est une constante qui ne d\'epend que de
$\alpha$ et de $\beta$.
\end{thm}
Ce th\'eor\`eme signifie que
$$\mathcal{C}_2(\delta,\aA^{+}_{\alpha,\beta})\leq\delta^{-c_{\alpha,\beta}},
\qquad0<\delta\leq1.$$
 Compte tenu du th\'eor\`eme \ref{symet2}, la
preuve du th\'eor\`eme \ref{2-d-ab} est une cons\'equence directe du
lemme suivant:

\begin{lem}\label{dnv}  Soient $\alpha,\beta>0.$  On a
$$\mathcal{C}_{2}(\delta,\aA^{+}_{\alpha,\beta})\leq c\
\cun^{3}(\delta^{2},\aA_{\alpha,\beta}),\qquad 0<\delta\leq1,$$ o\`u
$c$ est une constante universelle.
\end{lem}

 Pour chaque fonction
$g\in\aA(\D^2)$ et pour $0\leq r,s\leq1$  on d\'efinit la fonction $(g)_{r,s}\in\aA(\D^{2})$  par
$$(g)_{r,s}(z,w)=g(rz,sw)\qquad z,w\in{\D}. $$
Soit $\partial_{z}=\partial/\partial z$ la d\'eriv\'ee partielle par
rapport \`a la  variable $z$ et on d\'esigne par $\text{\rm Hol
}(\mathcal{U})$ l'ensemble des fonctions holomorphes sur l'ouvert
$\mathcal{U}$. Nous avons besoin du sous--lemme suivant:

\begin{slem}\label{lem1}
Soient $\mathcal{U}$ un voisinage ouvert de $\overline{\D}^2$ et
$a,b\in \mathcal{C}^{\infty}(\mathcal{U})$ tels que
$$b(z,w)=-\frac{1}{\pi}
\int\limits_{\mathbb{D}}\frac{a(\xi,w)}{\xi-z} dA(\xi),\qquad
z,w\in\D,
$$ o\`u $dA$ est la mesure d'aire. On a
$$\begin{array}{lll}
\widehat{b}(n,m)=0,&
(n,m)\in{\N}\times{\Z},\\
\displaystyle\widehat{b}(n,m)=2\int_{0}^{1}\widehat{(a)_{r,1}}(n+1,m)r^{-n}dr,&
(n,m)\not\in
{\N}\times{\Z}.\\
\end{array}
$$
De plus, si $f\in \text{\rm Hol }(\mathcal{U}),$ alors
\begin{align*}
\sum\limits_{n,m=0}^{\infty}|\widehat{fb}&(n,m)|(1+n)^{\alpha}(1+m)^{\beta}
\leq 2\| f\|_{\alpha,\beta} \int\limits_{0}^{1}\|
(a)_{r,1}\|_{\alpha,\beta}dr.
\end{align*}
\end{slem}

\begin{proof} Par le th\'eor\`eme de convergence domin\'ee on a
$$b(e^{i\theta},e^{i\varphi})=-\frac{1}{\pi}
\int\limits_{\mathbb{D}}\frac{a(\xi,e^{i\varphi})}{\xi-e^{i\theta}}
dA(\xi).$$
D'apr\`es le th\'eor\`eme de Fubini on obtient:
$$\widehat{b}(n,m) =
-\frac{1}{2\pi^{2}}\int\limits_{\mathbb{D}}\Big(\int\limits_{0}^{2\pi}
a(\xi,e^{i\varphi})e^{-im\varphi} d\varphi\Big)
\Big(\frac{1}{2\pi}\int\limits_{0}^{2\pi}
\frac{e^{-in\theta}}{\xi-e^{i\theta}}d\theta\Big)
d\mathrm{A}(\xi).$$
Si $(n,m)\in\mathbb{N}\times\mathbb{Z},$ alors
$\widehat{b}(n,m)=0.$ Sinon
\begin{multline*}\widehat{b}(n,m)
=\frac{1}{2\pi^2}\int\limits_{\mathbb{D}}
\Big(\int\limits_{0}^{2\pi}a(\xi,e^{i\varphi})e^{-im\varphi}d\varphi\Big)
\xi^{-n-1}d\mathrm{A}(\xi)\\ =
2\int\limits_{0}^{1}\Big(\frac{1}{4\pi^2} \int\limits_{0}^{2\pi}
\int\limits_{0}^{2\pi}
a(re^{i\theta},e^{i\varphi})e^{-i(n+1)\theta}e^{-im\varphi}d\theta
d\varphi\Big) r^{-n}dr \\=
2\int\limits_{0}^{1}\widehat{(a)_{r,1}}(n+1,m)r^{-n}dr.
\end{multline*}
Puisque
$$\widehat{fb}(n,m)=\sum\limits_{l=1}^{\infty}\sum\limits_{k=-\infty}^{m}
\widehat{b}(-l,k)\widehat{f}(n+l,m-k),$$ alors
\begin{multline}
 \sum\limits_{n,m=0}^{\infty}|\widehat{fb}(n,m)|(1+n)^{\alpha}(1+m)^{\beta}\\
=\sum\limits_{n,m=0}^{\infty}|
\sum\limits_{l=1}^{\infty}\sum\limits_{k\leq m}
\widehat{b}(-l,k)\widehat{f}(n+l,m-k)(1+n)^{\alpha}(1+m)^{\beta}|
\\
\leq  \sum\limits_{n,m=0}^{\infty}
\sum\limits_{l=1}^{\infty}\sum\limits_{k\leq m}|\widehat{b}(-l,k)|
(1+|k|)^\beta |\widehat{f}(n+l,m-k)|(1+n+l)^{\alpha}(1+m-k)^{\beta}\\
\leq\| f\|_{\alpha,\beta}\sum\limits_{l=1}^{\infty}
\sum\limits_{-\infty}^{+\infty}|\widehat{b}(-l,k)|
(1+|k|)^\beta\\
 \leq
\|f\|_{\alpha,\beta}\sum\limits_{l=1}^{\infty}\sum\limits_{-\infty}^{\infty}
\int\limits_{0}^{1}|\widehat{(a)_{r,1}}(1-l,k)|(1+|k|)^\beta
r^{l}dr  \\
\label{eq5}\leq 2\| f\|_{\alpha,\beta}
\int\limits_{0}^{1}\| (a)_{r,1}\|_{\alpha,\beta}dr.
\end{multline}
\end{proof}

{\textbf{Preuve du lemme \ref{dnv}}}. Soit $\delta\in ]0,1]$
et soient $f_{1}, f_{2}\in \aA_{\alpha,\beta}^{+}$ satisfaisant
\begin{equation}\label{eq2}\left\{
\begin{array}{l}
\|f_{1}\|^2_{\alpha,\beta}+\|f_{2}\|^2_{\alpha,\beta}\leq 1,  \\
\mathrm{F}(z,w)=|f_{1}(z,w)|^{2}+|f_{2}(z,w)|^{2}\geq\delta^{2}.
\end{array}
\right.
\end{equation}
Supposons d'abord que $f_1$ et $f_2$ sont holomorphes au voisinage
de $\overline{{\D}}^{2}.$ On pose
$$\phi_{i}=\frac{\overline{f}_{i}}{\mathrm{F}},\qquad i=1,2.$$
Comme dans la r\'esolution du probl\`eme de la couronne \cite{Ko}, nous allons
corriger la solution $(\phi_{1},\phi_{2})$ de l'\'equation
$f_{1}\phi_{1}+f_{2}\phi_{2}=1,$ pour aboutir a une solution
holomorphe au voisinage du bidisque. D'abord, on commence par
corriger les fonctions  $\phi_{1}$ et $\phi_{2}$ par rapport \`a la
premi\`ere variable $z$, soit
$$\left\{
\begin{array}{lll}
g_{1}&=&\phi_{1}+ f_{2}b,\\
g_{2}&=&\phi_{2}- f_{1}b,\\
 \overline{\partial}_{z}b&=&\phi_{1}\overline{\partial}_{z}\phi_{2}-
\phi_{2}\overline{\partial}_{z}\phi_{1}=a.
\end{array}
\right.
$$ Il est facile de
v\'erifier que  $f_{1}g_{1}+f_{2}g_{2}=1$ et que les fonctions
$g_{i}$ $(i=1,2)$ sont holomorphes sur le disque, par rapport \`a la
premi\`ere variable. Ensuite, nous corrigeons \`a nouveau les fonctions obtenues
$g_{1}$ et $g_{2}$ comme suit
$$\left\{
\begin{array}{lll}
h_{1}&= &g_{1} + f_{2}d,\\
 h_{2}&= &g_{2} -f_{1}d,\\
 \overline{\partial}_{w}d&=&g_{1}\overline{\partial}_{w}g_{2}-
g_{2}\overline{\partial}_{w}g_{1}=c.
\end{array}
\right.
$$
 On obtient ainsi  des fonctions holomorphes sur le bidisque
$h_{1},h_{2}$ satisfaisant l'identit\'e de Bezout suivante
$f_{1}h_{1} + f_{2}h_{2}=1.$ \\
Les solutions $b$ et $d$ sont donn\'ees par:
$$
b(z,w)=-\frac{1}{\pi}
\int\limits_{\mathbb{D}}\frac{a(\xi,w)}{\xi-z} dA(\xi)\quad \ \mbox{et}\quad \ d(z,w)=-\frac{1}{\pi}
\int\limits_{\mathbb{D}}\frac{c(z,\zeta)}{\zeta-w} dA(\zeta)\quad
.$$
Dans ce qui suit nous allons majorer les normes de $h_{1}$ et de
$h_{2}$. On a
\begin{eqnarray}\nonumber
\| h_{1}\|_{\alpha,\beta} &=&
\|\phi_{1}+ f_{2}b+ f_{2}d\|_{\alpha,\beta}\\\nonumber&\leq&
\|\phi_{1}
\|_{\alpha,\beta}+
\sum\limits_{n,m=0}^{\infty}|\widehat{f_{2}b}(n,m)|(1+n)^{\alpha}(1+m)^{\beta}
\\\label{eqq3}&&+\sum\limits_{n,m=0}^{\infty}|\widehat{f_{2}
d}(n,m)|(1+n)^{\alpha}(1+m)^{\beta},
\end{eqnarray}
On a
\begin{equation}\label{eqq4}
\| \phi_{1}\|_{\alpha,\beta}\leq
\|\mathrm{F}^{-1}\|_{\alpha,\beta}\|
f_{1}\|_{\alpha,\beta}\leq\cun(\delta^{2},\aA_{\alpha,\beta})\|
f_{1}\|_{\alpha,\beta}.
\end{equation}
D'apr\'es le sous--lemme \ref{lem1}, on obtient
\begin{eqnarray}\label{eqq5}\sum\limits_{n,m=0}^{\infty}|\widehat{f_{2}b}(n,m)|(1+n)^{\alpha}(1+m)^{\beta}
\leq 2\| f_{2}\|_{\alpha,\beta} \int\limits_{0}^{1}\|
(a)_{r,1}\|_{\alpha,\beta}dr.
\end{eqnarray}
Puisque
\begin{multline*}
\| (a)_{r,1}\|_{\alpha,\beta}= \| \left(
\overline{(\partial_{z}f_{1})_{r,1}(f_{2})_{r,1}}-
\overline{(f_{1})_{r,1}(\partial_{z}f_{2})_{r,1}}
\right)(\mathrm{F}^{-2})_{r,1}\|_{\alpha,\beta} \\
\leq \| (\mathrm{F}^{-1})_{r,1}\|_{\alpha,\beta}^{2} \left( \|
(\partial_{z}f_{1})_{r,1}\|_{\alpha,\beta} \|
(f_{2})_{r,1}\|_{\alpha,\beta}+ \|
(f_{1})_{r,1}\|_{\alpha,\beta} \|
(\partial_{z}f_{2})_{r,1}\|_{\alpha,\beta}
\right)\\
\leq \cun^2(\delta^{2},\aA_{\alpha,\beta})\left( \|
(\partial_{z}f_{1})_{r,1}\|_{\alpha,\beta} \| f_{2}\|_{\alpha,\beta}+\| f_{1}
\|_{\alpha,\beta}\| (\partial_{z}f_{2})_{r,1} \|_{\alpha,\beta} \right),
\end{multline*}
donc
\begin{multline*}
\int\limits_{0}^{1}\| (a)_{r,1}\|_{\alpha,\beta}dr\leq
\cun^2(\delta^{2},\aA_{\alpha,\beta})\times\\
\times\left( \| f_{2}\|_{\alpha,\beta}\int\limits_{0}^{1}\|
(\partial_{z}f_{1})_{r,1}\|_{\alpha,\beta}dr+\|
f_{1}\|_{\alpha,\beta}\int\limits_{0}^{1}\|
(\partial_{z}f_{2})_{r,1}\|_{\alpha,\beta} dr\right).
\end{multline*}
Pour $i=1,2$
\begin{eqnarray*}
\int\limits_{0}^{1}\|
(\partial_{z}f_{i})_{r,1}\|_{\alpha,\beta}dr&=&\int\limits_{0}^{1}
\sum\limits_{n,m=0}^{\infty}n|\widehat{f_{i}}(n,m)|(1+n)^\alpha(1+m)^{\beta}
r^{n-1}dr\\&\leq&\| f_{i}\|_{\alpha,\beta}.
\end{eqnarray*}
Ce qui donne
\begin{equation} \label{eq6}\int\limits_{0}^{1}\|
(a)_{r,1}\|_{\alpha,\beta}dr\leq 2\
\cun^2(\delta^{2},\aA_{\alpha,\beta})\| f_{1} \|_{\alpha,\beta}\|
f_{2}\|_{\alpha,\beta}.
\end{equation}
En combinant les in\'egalit\'es \eqref{eqq5} et \eqref{eq6}, on obtient
\begin{eqnarray}\label{eq7}
\sum\limits_{n,m=0}^{\infty}|\widehat{f_{2}b}(n,m)|(1+n)^{\alpha}(1+m)^{\beta}\leq4\
\cun^2(\delta^{2},\aA_{\alpha,\beta})\| f_{1}\|_{\alpha,\beta}.
\end{eqnarray}
Dans ce qui suit nous allons majorer la quantit\'e suivante
$$\sum\limits_{n,m=0}^{\infty}|\widehat{f_{2}
d}(n,m)|(1+n)^{\alpha}(1+m)^{\beta}.$$
En utilisant le sous--lemme
\ref{lem1} on obtient
\begin{eqnarray} \label{eq8}\sum\limits_{n,m=0}^{\infty}|\widehat{f_{2}d}(n,m)
|(1+n)^{\alpha}(1+m)^{\beta} \leq  2\|f_{2}\|_{\alpha,\beta}
\int\limits_{0}^{1}\| (c)_{1,s}\|_{\alpha,\beta}ds,
\end{eqnarray}
o\`u
$$c=g_{1}\overline{\partial_{w}}g_{2}-g_{2}\overline{\partial_{w}}g_{1}
=\phi_{1}\overline{\partial_{w}}\phi_{2}-\phi_{2}\overline{\partial_{w}}\phi_{1}-
\overline{\partial_{w}}b.$$
De la m\^eme fa\c{c}on que dans l'\'equation
\eqref{eq6} on peut montrer que
\begin{equation}\label{eq9}\int\limits_{0}^{1}\|\Big(\phi_{1}\overline{\partial_{w}}\phi_{2}-\phi_{2}
\overline{\partial_{w}}\phi_{1}
\Big)_{1,s}\|_{\alpha,\beta}ds\leq 2\
\cun^2(\delta^{2},\aA_{\alpha,\beta})\| f_{1} \|_{\alpha,\beta}\|
f_{2}\|_{\alpha,\beta}.
\end{equation}
Puisque
$$\overline{\partial _w}b(z,sw)=-\frac{1}{\pi}
\int\limits_{\mathbb{D}}\frac{\overline{\partial _w}a(\xi,sw)}{\xi-z} dA(\xi),\qquad
z,w\in\D,$$ d'apr\`es le sous--lemme \ref{lem1}, on obtient
$$
\begin{array}{lll}
\widehat{(\overline{\partial_{w}}b)_{1,s}}(n,m)=0,&(n,m)\in{\N}\times{\Z},\\
\displaystyle\widehat{(\overline{\partial_{w}}b)_{1,s}}(n,m)
\leq 2\displaystyle\int_{0}^{1}|\widehat{(\overline{\partial_{w}}(a)_{r,s}}|(n+1,m)r^{-n}dr,&
(n,m)\not\in
{\N}\times{\Z}.\\
\end{array}
$$
Ce qui entra\^{\i}ne que
$$\int\limits_{0}^{1}\|
(\overline{\partial_{w}}b)_{1,s}\|_{\alpha,\beta}ds\leq
2\int\limits_{0}^{1}\int\limits_{0}^{1}\|
(\overline{\partial_{w}}a)_{r,s}\|_{\alpha,\beta}drds.$$ Comme
$a=\Big(\overline{\partial_{z}(f_{1})f_{2}}-\overline{f_{1}
\partial_{z}f_{2}}\Big)\mathrm{F}^{-2}$,
on a
\begin{multline*}
\overline{\partial_{w}}a=
\Big(\overline{\partial_{z}\partial_{w}(f_{1})f_{2}+\partial_{z}
f_{1}\partial_{w}f_{2}
-\partial_{w}f_{1}\partial_{z}f_{2}-f_{1}\partial_{z}
\partial_{w}f_{2}}\Big)\mathrm{F}^{-2}-\\
2\Big(\overline{\partial_{z}(f_{1})f_{2}}-\overline{f_{1}
\partial_{z}(f_{2})}\Big)\Big(f_{1}\overline{\partial_{w}f_{1}}+
f_{2}\overline{\partial_{w}f_{2}}\Big)\mathrm{F}^{-3}.
\end{multline*}
Or,
\begin{multline*}
\|
\Big(\overline{\partial_{z}\partial_{w}(f_{1})f_{2}+\partial_{z}f_{1}
\partial_{w}f_{2}
-\partial_{w}f_{1}\partial_{z}f_{2}-f_{1}\partial_{z}
\partial_{w}f_{2}}\Big)_{r,s}\Big(\mathrm{F}^{-2}\Big)_{r,s}\|_{\alpha,\beta}
\\\leq
\cun^2(\delta^{2},\aA_{\alpha,\beta})\Big(\|\partial_{z}
\partial_{w}(f_{1})_{r,s}\|_{\alpha,\beta}\|
f_{2}\|_{\alpha,\beta}+\|\partial_{z}(f_{1})_{r,1}
\|_{\alpha,\beta}\|\partial_{w}(f_{2})_{1,s} \|_{\alpha,\beta}
\\+
\|\partial_{w}(f_{1})_{1,s}\|_{\alpha,\beta}\|\partial_{z}(f_{2})_{r,1}
\|_{\alpha,\beta}+ \|
f_{1}\|_{\alpha,\beta}\|\partial_{z}\partial_{w}(f_{2})_{r,s}
\|_{\alpha,\beta}\Big).
\end{multline*} De plus,
\begin{multline*}
\|(\overline{\partial_{z}(f_{1})f_{2}}-\overline{f_{1}
\partial_{z}(f_{2})})_{r,s}(f_{1}\overline{\partial_{w}f_{1}}+
f_{2}\overline{\partial_{w}f_{2}})_{r,s}(\mathrm{F}^{-3})_{r,s}\|_{\alpha,\beta}
\\\leq
\cun^3(\delta^{2},\aA_{\alpha,\beta})\Big(\|(\partial_{z}(f_{1}))_{r,1}
\|_{\alpha,\beta}\| f_{2}\|_{\alpha,\beta}+\|
f_{1}\|_{\alpha,\beta}\|\partial_{z}((f_{2}))_{r,1}\|_{\alpha,\beta}\Big)\times
\\\times\Big(\|
f_{1}\|_{\alpha,\beta}\|(\partial_{w}f_{1})_{1,s} \|_{\alpha,\beta}+ \|
f_{2}\|_{\alpha,\beta}\|(\partial_{w}f_{2})_{1,s} \|_{\alpha,\beta}\Big).
\end{multline*}
Ce qui donne
\begin{multline*}
\int\limits_{0}^{1}\int\limits_{0}^{1}\|
(\overline{\partial_{w}}a)_{r,s}\|_{\alpha,\beta}drds\\
\leq 4\
\cun^2(\delta^{2},\aA_{\alpha,\beta})\| f_{1}\|_{\alpha,\beta}\| f_{2}
\|_{\alpha,\beta}
 + 2\ \cun^3(\delta^{2},\aA_{\alpha,\beta})\| f_{1}\|_{\alpha,\beta}\|
f_{2}\|_{\alpha,\beta}\\
\leq 6\ \cun^3(\delta^{2},\aA_{\alpha,\beta})\|
f_{1}\|_{\alpha,\beta}\| f_{2}\|_{\alpha,\beta}.
\end{multline*} Donc
\begin{equation}\label{eq10} \int\limits_{0}^{1}\|
(\overline{\partial_{w}}b)_{1,s}\|_{\alpha,\beta}ds\leq12\
\cun^3(\delta^{2},\aA_{\alpha,\beta})\| f_{1}\|_{\alpha,\beta}\|
f_{2}\|_{\alpha,\beta}.
\end{equation}
En combinant \eqref{eq8}, \eqref{eq9} et \eqref{eq10}, on obtient
\begin{equation}\label{eq11}
\sum\limits_{n,m=0}^{\infty}|\widehat{f_{2}d}(n,m)|(1+n)^{\alpha}(1+m)^{\beta}\leq
28\ \cun^3(\delta^{2},\aA_{\alpha,\beta})\| f_{1}\|_{\alpha,\beta}.
\end{equation}
Les in\'egalit\'es \eqref{eqq3}, \eqref{eqq4}, \eqref{eq7} et
\eqref{eq11} entra\^{\i}nent que
$$
\| h_{1}\|_{\alpha,\beta} \leq 33\
\cun^3(\delta^{2},\aA_{\alpha,\beta})\| f_{1}\|_{\alpha,\beta} .
$$
De la m\^eme fa\c{c}on on peut montrer que
$$\| h_{2}\|_{\alpha,\beta}
\leq 33\ \cun^3(\delta^{2},\aA_{\alpha,\beta}) \|
f_{2}\|_{\alpha,\beta}.
$$
On obtient donc
$$(\| h_{1}\|^{2}_{\alpha,\beta}+\|
h_{2}\|^{2}_{\alpha,\beta})^{1/2}\leq 33\
\cun^3(\delta^{2},\aA_{\alpha,\beta}).$$

On suppose maintenant que  $f_1$ et $f_2$ ne sont pas holomorphes au
voisinage de $\overline{{\D}}^{2}.$ On consid\`ere les fonctions
$(f_1)_{r,r}$ et $(f_2)_{r,r}$ $(r<1).$ Donc il existe
 $h_{1,r},h_{2,r}\in\aA_{\alpha,\beta}^{+}$ tel que
$$\left\{\begin{array}{l}
 (f_1)_{r,r}h_{1,r}+(f_2)_{r,r}h_{2,r}=1,\\\\
(\|h_{1,r}\|^{2}_{\alpha,\beta}+\|h_{2,r}\|^{2}_{\alpha,\beta})^{1/2}\leq
33\ \cun^3(\delta^{2},\aA_{\alpha,\beta}),\quad r<1.
\end{array}
\right.
$$
Puisque $\aA_{\alpha,\beta}^{+}$ peut \^etre identifi\'e comme
 le dual de
$$\ell^{\infty}_{\alpha,\beta}(\N^2):=\{(u_{n,m})_{n,m\geq0}\text{ :
}\sup\limits_{n,m\geq
0}\frac{|u_{n,m}|}{(1+n)^{\alpha}(1+m)^{\beta}}<\infty\}.$$
La boule
unit\'e ferm\'ee de $\aA_{\alpha,\beta}^{+}$ est compacte pour
la topologie $\ast-$faible. On en d\'eduit qu'il existe
$h_1,h_2\in\aA_{\alpha,\beta}^{+}$ tels que
$$ \left\{\begin{array}{l}
f_1h_1+f_2h_2=1,\\\\
(\| h_{1}\|^{2}_{\alpha,\beta}+\|
h_{2}\|^{2}_{\alpha,\beta})^{1/2}\leq 33\
\cun^3(\delta^{2},\aA_{\alpha,\beta}).\\
\end{array}
\right.
$$ Ce qui termine la preuve du
lemme \ref{dnv}.


\section{ Preuve du Theoreme 1}

\begin{lem}\label{lem15}
Soit $f\in \aA(\D^2)$ ne s'annulant pas sur $\D\times\overline{\D}.$
Alors il existe $M>0$ tel que
$$\frac1{|f(z, w)|}\leq e^{\frac{M}{1-| z |}},\qquad z, w \in\D.$$
\end{lem}

\begin{proof}
Supposons que $\|f\|_\infty\leq 1$. 
 Puisque l'application $ \log |f(\cdot ,w)|^{-1}$ est positive et harmonique sur $\D$,
il existe une mesure positive $\mu=\mu _w$ sur $\T$ telle que
$$\log|f(z,w)|^{-1}=\frac{1}{2\pi}\int_{-\pi}^{\pi}\frac{1-|z|^2}{|e^{i\varphi}-z|^2}
\log|f(e^{i\varphi} ,w)|^{-1}d\mu (\varphi)\qquad(w\in\D).$$ On
obtient donc
$$\log|f(z,w)|^{-1}\leq 2\log|f(0,w)|^{-1}(1-|z|)^{-1}.$$
Par cons\'equent
$$\log|f(z,w)|^{-1}\leq M(1-|z|)^{-1},$$
 o\`u  $M=2\sup\{\log|f(0,w)|^{-1}\ :\ w\in \overline{\D} \}.$
 \end{proof}
Le r\'esultat suivant est d\^u \`a H.Hedenmalm
\cite{Hed1}, proposition 1.2.
\begin{lem}\label{prop1}
Soit $f\in\aA(\mathbb{D}^{2})$ telle que
$Z_f\subset\{1\}\times\overline{\D}$. Alors les trois assertions
suivantes sont \'equivalentes:
\begin{enumerate}
\item $f(\cdot,w)$ est une fonction ext\'erieure pour tout $w\in\overline{{\D}},$
\item $\displaystyle\lim_{r\to 1^{-}}(1-r)\log|f(r,w)|=0$ pour tout
$w\in\overline{{\D}},$
\item $\displaystyle\lim_{r\to 1^{-}}\inf_{w\in
\overline{{\D}}}\ (1-r)\log|f(r,w)|=0.$
\end{enumerate}
\end{lem}

Nous aurons besoin du principe de Phragm\'en--Lindel\"of suivant (voir \cite{Hed3}).
\begin{lem}\label{dtw}
Soit $\phi$ une fonction analytique sur ${\C}\setminus\{1\}$ et soit
$\varepsilon>0.$ On suppose que
$$
\begin{array}{llllr}
(i)\qquad\displaystyle &|\phi(\lambda)|\leq \displaystyle\frac{c_1}{(|\lambda|-1)^{N}}&,& \qquad1<|\lambda|<2,\\
(ii)\qquad\displaystyle &  |\phi(\lambda)| \leq c_2\displaystyle\exp\frac{c_3}{1-|\lambda|}&, &\qquad |\lambda|<1,\\
(iii)\qquad\displaystyle & |\phi(x)|\leq c_{\varepsilon}\displaystyle\exp\frac{\varepsilon}{1-x}&,& \qquad x<1,
\end{array}
$$
o\`u les $c_1,c_2,c_2$  sont des constantes positives, $N$ est
un entier et $c_{\varepsilon}$ est une constante qui  d\'epend de $\varepsilon.$ Alors $\phi$ est un polyn\^ome en
$\frac{1}{1-\lambda}$ de degr\'e inf\'erieur ou \'egal \`a $N.$
\end{lem}

{\bf Preuve du Th\'eor\`eme 1.} Supposons que
$I_{\aA^{+}_{\alpha,\beta}}(f)=I_{\aA^{+}_{\alpha,\beta}}(Z_f)$. Il
est clair que pour tout $w\in \overline{\D}$ on a
$I_{\aA(\D)}(f(\cdot,w))=I_{\aA(\D)}(Z_{f(\cdot,w)})$ et donc $f(\cdot, w)$ est ext\'erieure.

Supposons maintenant que $\|f\|_{{\alpha,\beta}}\leq 1$ et que  $f(\cdot,w)$ est  ext\'erieure pour tout $w\in\overline{\D}$.  Notons par $\pi$ la surjection canonique sur
$\aA^{+}_{\alpha,\beta}$ associ\'ee \`a
$I_{\aA^{+}_{\alpha,\beta}}(f)$, et par $u$ la fonction d\'efinie
par $u(z,w)=z$. Comme $Z_f = \{1\}\times \overline {\D}$, le spectre
de $\pi(u)$ est r\'eduit au singleton $\{1\}$. Par cons\'equent
l'application
$$\phi(\lambda)=(\lambda-\pi(u))^{-1}$$
est analytique sur $\C\setminus\{1\}.$ Nous allons montrer que
$\phi$ satisfait les trois conditions du Lemme \ref{dtw}.
Pour $1<|\lambda|<2$ on a:
\begin{multline}
\|\phi(\lambda)\|\leq\sum\limits_{n\geq0}\|
\pi(u)^{n}\lambda^{-n-1}\|_{{\alpha,\beta}}
\leq \sum\limits_{n\geq0}| \lambda|^{-n-1}(1+n)^{\alpha}
\label{eq1} \leq \frac{c}{(|\lambda|-1)^{1+\alpha}}.
\end{multline}
Pour $|\lambda|<1$, on pose
\begin{equation}\label{eq2}
 L_{\lambda}(f)(z,w)= \left \{ \begin{array}{ll}
\displaystyle \frac{f(z,w)-f(\lambda,w)}{z-\lambda}\quad si \ z\neq\lambda\\
\displaystyle\frac{\partial}{\partial z} f(\lambda,w) \qquad \qquad
si\ z=\lambda.
\end{array} \right.
\end{equation}
On obtient
 \begin{multline*} \| L_{\lambda}(f)\|_{\alpha,\beta}=\|
\sum\limits_{n,m=0}^{\infty}a_{n,m}(\frac{z^{n}-\lambda^{n}}{z-\lambda})w^{m}\|_{{\alpha,\beta}}=\\
\|\sum\limits_{n,m=0}^{\infty}a_{n,m}(z^{n-1}+
z^{n-2}\lambda+\cdots+\lambda^{n-1})w^{m}\|_{{\alpha,\beta}}
\leq\frac{\| f\|_{{\alpha,\beta}}}{1-|\lambda|} \leq\frac{1}{1-|
\lambda|}.
\end{multline*}
De l'\'equation \eqref{eq2}, on a
$$\|\phi(\lambda)\|\leq
\| \pi(L_{\lambda}(f))\|_{{\alpha,\beta}}\|
(f(\lambda,\cdot))^{-1}\|_{\beta}\leq \frac{1}{1-|\lambda|}\|
(f(\lambda,\cdot))^{-1}\|_{\beta}.$$
On pose
$$\delta(\lambda)=\inf\{|f(\lambda,w)|: \quad\ w\in\mathbb{D}
\}.$$
Puisque $\cun(\delta,\aA_\beta)\leq c/\delta^{2+1/\beta}$ on obtient
\begin{equation}\label{eq30}\|\phi(\lambda)\|\leq
\frac{C}{(1-|\lambda|)(\delta(\lambda))^{2+ 1/\beta}}.
\end{equation} D'apr\'es le lemme \ref{lem15}, il existe une constante $M>0$ telle que
\begin{equation}\label{eq+}\|\phi(\lambda)\|\leq
e^{\frac{M}{1-|\lambda|}}, \qquad |\lambda|<1.
\end{equation}
 Puisque $f(\cdot,w)$ est ext\'erieure pour tout $w\in \overline{\D}$, d'apr\'es le lemme \ref{prop1}, pour tout $\varepsilon>0$ , il
existe une constante $c_{\varepsilon}>0$ telle que
\begin{equation}\label{eq4}\delta(r)\geq
c_{\varepsilon}e^{\frac{-\varepsilon}{1-r}},  \qquad 0<r<1.
\end{equation}  On obtient alors
\begin{equation}\label{iii}\|\phi(r)\|\leq
c_{\varepsilon}e^{\frac{\varepsilon}{1-r}}, \qquad 0<r<1.
\end{equation}Soit
$g\in  (\aA^{+}_{\alpha,\beta}/I_{\aA^{+}_{\alpha,\beta}}(f))^{*},$ le dual de
$\aA^{+}_{\alpha,\beta}/I_{\aA^{+}_{\alpha,\beta}}(f)$. D'apr\'es les \'equations \eqref{eq1},
\eqref{eq+} et \eqref{iii}, la fonction
$$\lambda\longrightarrow
\phi_{g}(\lambda):=\langle\phi(\lambda),g\rangle$$ v\'erifie les
conditions du lemme \ref{dtw}. On en d\'eduit que $\phi_{g}$ est un
polyn\^ome en $\frac{1}{1-\lambda}$ de degr\'e $1.$ Donc
$\pi(1-u)=0$ et $I_{\aA^{+}_{\alpha,\beta}}(f)\supseteq
I_{\aA^{+}_{\alpha,\beta}}(\{1\}\times\overline{\mathbb{D}})$.
L'autre inclusion est triviale. Ce qui termine la preuve du
th\'eor\`eme \ref{theo11}.

\section{ Preuve du Theoreme 2}
Commen\c{c}ons d'abord par \'etablir quelque lemmes qui nous serons utile
dans la preuve du th\'eor\`eme \ref{theo12}.

\begin{lem} \label{apartien} Soit $f\in{\rm H}^{\infty}(\D^{2})$ et soient $\alpha,\beta\in [0,1[$. On suppose que
$f(\cdot,0)$, ${\partial_z}f(\cdot,0)\in \aA_\alpha^+$ et que
$f(0,\cdot)$, ${\partial_w}f(0,\cdot)\in \aA_\beta^+$. S'il existe
 $0<\varepsilon<2\inf{\{(1-\alpha),(1-\beta)\}}$ tel que
$$  \int_{\D^2}|\frac{\partial^4}{\partial z^2\partial w^2}f(z,w)|^{2}(1-|z|)^{2(1-\alpha)-\varepsilon}
  (1-|w|)^{2(1-\beta)-\varepsilon}dA(z)dA(w)<+\infty,$$
  alors  $f\in\aA^{+}_{\alpha,\beta}.$
\end{lem}

\begin{proof}
Nous avons
\begin{multline*}
f(z,w)=\sum_n\Big(\sum_m a_{n,m}w^m\Big)z^n\\=
\sum_{n\geq2}\Big(\sum_{m\geq2}a_{n,m}w^m\Big)z^n
+\sum_{n\geq0}a_{n,0}z^n+\sum_{n\geq0}a_{n,1}z^nw+\\
+\sum_{m\geq2}a_{0,m}w^m +\sum_{m\geq2}a_{1,m}zw^m.\\=
\sum_{n\geq2}\Big(b_n(w)\Big)z^n + f(z,0)+ w{\partial_w}f(z,0)
\\+(f(0,w)-a_{0,0} -a_{0,1}w)+z{\partial_z}f(0,w)-a_{1,0} -a_{1,1}w).
\end{multline*}
o\`u $b_n(w)=\sum\limits_{m\geq2}a_{n,m}w^m.$  Donc 
\begin{multline*}
   \sum_{n\geq2}\sum_{m\geq2}|a_{n,m}|(1+n)^{\alpha}(1+m)^{\beta}=
   \sum_{n\geq2}\Big(\sum_{m\geq2}|a_{n,m}|(1+m)^{\beta}\Big)(1+n)^{\alpha}\\
\leq
   \Big(\sum_{m\geq2}\frac{1}{(1+m)^{1+\epsilon}}\Big)^{1/2}\sum_{n\geq2}\Big(\sum_{m\geq2}|a_{n,m}|^2(1+m)^{2\beta+1+\epsilon}\Big)^{1/2}(1+n)^{\alpha}\\
\leq c_{\epsilon,\beta}
\sum_{n\geq2}\Big(\int_{\D}|\frac{\partial^2 b_{n}}{\partial w^2}(w)|^{2}(1-|w|^{2})^{2(1-\beta)-\epsilon}dA(w)\Big)^{1/2}(1+n)^{\alpha} \\
\leq c_{\epsilon,\beta}
\Big(\sum_{n\geq2}\frac{1}{(1+n)^{1+\epsilon}}\Big)^{1/2}\times\\
\times\Big(\int_{\D}\sum_{n\geq2}\Big(|\frac{\partial^2 b_{n}}{\partial w^2}(w)|^2(1+n)^{2\alpha+1+\epsilon}\Big)
(1-|w|^{2})^{2(1-\beta)-\epsilon}dA(w)\Big)^{1/2}\\
\leq c_{\epsilon,\alpha, \beta}\int_{\D^2}|\frac{\partial^4
f}{\partial z^2\partial
w^2}(z,w)|^{2}(1-|z|^{2})^{2(1-\alpha)-\epsilon}(1-|w|^{2})^{2(1-\beta)-\epsilon}
dA(z)dA(w).
\end{multline*}
Ce qui termine la preuve.

\end{proof}

\begin{lem}\label{lafu}
Soit $u$ la fonction d\'efinie sur $\D^2$ par
 $$u(z,w):=\frac{(1-z)(1-w)}{[(1-z)^{1/2}+(1-w)^{1/2}]^2}.$$
Alors $u$ v\'erifie les propri\'et\'es suivantes:
\begin{enumerate}
\item $|u(z,w)|\asymp \inf\{|1-z|,|1-w|\}$,
\item $\text{\rm Im }(u):=\{u(z,w)\text{ : }(z,w)\in\D^2\}\subset\{z\in\D\text{ : } |1-z|\leq1\}$,
\item $u^2\in \aA^{+}_{\alpha,\beta}$.
\end{enumerate}
\end{lem}
\begin{proof} Pour montrer (1), il suffit de remarquer que
\begin{multline*}
|(1-z)^{1/2}+(1-w)^{1/2}|\geq \text{\rm Re }[(1-z)^{1/2}+(1-w)^{1/2}]\\\geq \max\{\text{\rm Re }(1-z)^{1/2},
\text{\rm Re }(1-w)^{1/2}\}\geq
\cos \frac{\pi}{4}\max\{|1-z|^{1/2}, |1-w|^{1/2}\}.
\end{multline*}

 (2). Puisque $\{z\in\D\text{ : } |1-z|\leq1\}=\{z\in\C\text{ : }\text{\rm Re } (z)\geq|z|^2/2\}$,  pour $(z,w)\in\D^2$, $\text{\rm Re }(1/(1-z))\geq 1/2$ et $\text{\rm Re }(1/(1-w))\geq 1/2$ et
$$\text{\rm Re }\frac{u(z,w)}{|u(z,w)|^2}=\text{\rm Re }\Big[\Big(\frac{1}{1-z}\Big)^{1/2}+
\Big(\frac{1}{1-w}\Big)^{1/2}\Big]^2\geq\frac{1}{2}.$$
Donc $u(z,w)\in   \{z\in\D\text{ : } |1-z|\leq1\}$.

(3).  Nous avons
$$\frac{\partial^4}{\partial z^2\partial w^2}u^2(z,w)=\frac{105}{2}
\frac{(1-z)(1-w)}{[(1-z)^{1/2}+(1-w)^{1/2}]^8}-\frac{45}{2}\frac{(1-z)(1-w)}{[(1-z)^{1/2}+(1-w)^{1/2}]^6}.$$
On a
\begin{multline*}\int_{\D^2}|\frac{\partial^4}{\partial z^2\partial w^2}u^2(z,w)|^{2}
(1-|z|)^{1-\alpha}(1-|w|)^{1-\beta}dA(z)dA(w)\leq \\
c\int_{\D^2}\frac{1}{|1- z|^2}\frac{1}{|1-w|^2}(1-|z|)^{1-\alpha}(1-|w|)^{1-\beta}
dA(z)dA(w)\leq c_{\alpha,\beta}.
\end{multline*}
Le r\'esultat d\'ecoule alors du lemme \ref{apartien} avec
$\epsilon=\min \{1-\alpha, 1-\beta \}$.

\end{proof}

{\it \bf Preuve du th\'eor\`eme \ref{theo12} :}  Soit $u$ la fonction du  lemme \ref{lafu}. On pose $v:=u^2$, on a
$v\in\aA^{+}_{\alpha,\beta}.$ Soit $\pi: \aA^{+}_{\alpha,\beta}\to \aA^{+}_{\alpha,\beta}\big/I_{\aA^{+}_{\alpha,\beta}}(f)$
 la surjection canonique.  Le spectre de
$\pi(v)$ est r\'eduit \`a $\{0\}$, donc l'application
$$\varphi(\lambda)=(\lambda-\pi(v))^{-1}$$ est
analytique sur ${\C}\backslash\{0\}.$  Le lemme
\ref{lafu} entra\^{\i}ne que
$$\text{\rm Im }(v)\subset \cS:=\{z\in\C\setminus]-\infty,\
0[\text{ : } |z^{1/2}-1|\leq1\}.$$ D'apr\'es le Lemme
\ref{symet2} on a $\cun(\delta,\aA_{\alpha,\beta})\leq c\delta^{-c_{\alpha, \beta}}$.
 Donc il existe un entier
$N$  tel que pour tout $\lambda<0$ on a
\begin{equation}\label{est1}\| \varphi(\lambda)\| \leq \|(\lambda - v)^{-1}\|_{\aA^+ _{\alpha,\beta}}
\leq \frac{1}{ \text{\rm dist}(\lambda,\cS)^{N}} \leq \frac{c}{|\lambda|^{2N}},
\end{equation}
Soit $\lambda\neq 0.$ On pose
$$\delta(\lambda):=\inf \{|\lambda-v(z,w)|+|f(z,w)|\text{ : }(z,w)\in\D^2\}.$$
 Soit 
$$\cQ_{\lambda}:=\{(z,w)\in\overline{{\D}}^{2}:\
 |v(z,w)-\lambda|\leq{|\lambda|}/{2}\}.$$
Nous avons
$$|\lambda-v(z,w)|+|f(z,w)|\geq|\lambda-v(z,w)|\geq\frac{|\lambda|}{2},
\qquad (z,w)\notin\cQ_{\lambda}.$$
D'apr\'es le lemme \ref{lafu}, on a
$$ \inf\{|1-z|^2,|1-w|^2\}\geq c_1 |v(z,w)|\geq c_1|\lambda|/2,
\qquad (z,w)\in\cQ_{\lambda},$$
donc, pour tout  $\varepsilon>0$
$$|\lambda-v(z,w)| +|
 f(z,w)|\geq|f(z,w)|\geq
c_\varepsilon e^{-\varepsilon|\lambda|^{-1/2}}, \qquad
(z,w)\in\cQ_{\lambda}.$$
Par cons\'equent
\begin{equation}
\label{eqf}
\delta(\lambda)\geq
c_\varepsilon e^{-\varepsilon|\lambda|^{-1/2}}.
\end{equation}
D'apr\`es le
Lemme \ref{2-d-ab}, il existe $h_1,h_2\in \aA^{+}_{\alpha,\beta}$
tel que
 $$\left\{
 \begin{array}{l} (\lambda-v(z))h_1(z)+f(z)h_2(z)=1,\\ \\
 \|h_1\|^{2}_{\alpha,\beta}+\|h_2\|^{2}_{\alpha,\beta}\leq \delta(\lambda)^{-c_{\alpha,\beta}}.\\
 \end{array}
 \right.
 $$
  De l'in\'egalit\'e \eqref{eqf}, nous obtenons
\begin{multline}
\|\varphi(\lambda)\|_{\alpha,\beta}=\|(\lambda-\pi(v))^{-1}\|_{\alpha,\beta}=\|\pi(h_1)\|_{\alpha,\beta}\\
\leq \delta(\lambda)^{-c_{\alpha,\beta}/2}
\label{est2}\leq c_\varepsilon
e^{\varepsilon|\lambda|^{-1/2}}, \qquad \lambda\neq 0.
\end{multline}
Soit
$g\in(\aA^{+}_{\alpha,\beta}/I_{\aA^{+}_{\alpha,\beta}}(f))^{*}.$ On
d\'efinit $\varphi_g$ par
$$\lambda\longrightarrow
\varphi_{g}(\lambda):=<\varphi(\lambda),g>.$$
Le th\'eor\`eme
classique de Phragm\'en-Lindel\"of appliqu\'e \`a la fonction
$\varphi_{g}$ nous permet de d\'eduire, en utilisant les
in\'egalit\'es \eqref{est1} et \eqref{est2}, que $\varphi_{g}$ est
un polyn\^ome en $1/\lambda$ de degr\'e inf\'erieur \`a $2N.$  Le
th\'eor\`eme de Banach--Steinhaus nous donne  $\pi(v)^{2N}=0$ et
$v^{2N}\in I_{\aA^{+}_{\alpha,\beta}}(f).$ Soit maintenant la fonction $g$ d\'efinie par
$$g(z,w)= (1-z)(1-w)\big((1-z)^{1/2}+(1-w)^{1/2}\big)^{8N} \qquad
(z,w)\in{\D}^{2}.$$
Puisque  $$v^{2N}g(z,w)=(1-z)^{4N+1}(1-w)^{4N+1},$$
$(1-z)^{4N+1}(1-w)^{4N+1}\in I_{\aA^{+}_{\alpha,\beta}}(f)$. Pour conclure il suffit de remarquer que   pour tout $n$
 $$I_{\aA^{+}_{\alpha,\beta}}((1-z)^n(1-w)^n)=I_{\aA^{+}_{\alpha,\beta}}((1-z)(1-w))=
 I_{\aA^{+}_{\alpha,\beta}}(\{1\}\times\overline{\D}\cup\overline{\D}\times\{1\})$$ ce qui termine la preuve.

\subsection*{Remerciement} Les auteurs tiennent \`a remercier le professeur 
Alexander Borichev pour les discussions concernant ce travail.

\end{document}